\documentclass{amsart}

\usepackage[latin1]{inputenc}
\usepackage{latexsym}
\usepackage{amsfonts}
\usepackage{amssymb}
\usepackage{color}
\usepackage{dsfont}
\usepackage{varioref}
\usepackage{epsfig}
\usepackage{xypic}
\usepackage{hyperref}

\newtheorem{Lem}{Lemma}
\newtheorem{Exa}{Example}
\newtheorem{Def}[Lem]{Definition}
\newtheorem{Prop}[Lem]{Proposition}
\newtheorem{PropDef}[Lem]{Proposition/Definition}
\newtheorem{ThmDef}[Lem]{Theorem/Definition}

\newtheorem{Cor}[Lem]{Corollary}

\newcommand{\myinput}[1]{}

\newcommand{\dN}{{\mathds N}}
\newcommand{\dP}{{\mathds P}}
\newcommand{\dC}{{\mathds C}}
\newcommand{\dR}{{\mathds R}}

\newcommand{\smcdot}{{\textup{$\cdot$}}}
\newcommand{\ublb}[2]{\begin{tabular}{@{}c@{}}{\rule{0em}{1.6em}\setlength{\unitlength}{4144sp}%
\begin{picture}(347,215)(-41,616)
\thinlines
{\color[rgb]{0,0,0}\put(  1,659){\line( 2, 1){270}}
}%
\put(122,616){\makebox(0,0)[lb]{\smash{{\SetFigFont{7}{8.4}{\rmdefault}{\mddefault}{\updefault}{\color[rgb]{0,0,0}#2}%
}}}}
\put(143,756){\makebox(0,0)[rb]{\smash{{\SetFigFont{7}{8.4}{\rmdefault}{\mddefault}{\updefault}{\color[rgb]{0,0,0}#1}%
}}}}
\end{picture}%
}\end{tabular}}

\begin{document}

\title[Dessins d'Enfants and Hypersurfaces with Many $A_j$-Singularities]{Dessins
  d'Enfants and Hypersurfaces with\\Many $A_j$-Singularities
} 
\date{May 2005}
\author{Oliver Labs}
\address{Johannes Gutenberg Universit\"at Mainz, Germany}
\email{Labs@Mathematik.Uni-Mainz.de, mail@OliverLabs.net}

\subjclass{Primary 14J17, 14Q10}
\keywords{many cusps, many singularities, many nodes, dessins d'enfants}

\begin{abstract}
  We show the existence of surfaces of degree $d$ in
  $\dP^3(\dC)$ with approximately ${3j+2\over 6j(j+1)} d^3$ singularities of
  type $A_j, 2\le j\le d-1$.  
  The result is based on Chmutov's construction of nodal surfaces.
  For the proof we use plane trees related to the theory of Dessins d'Enfants.  

  Our examples improve the previously known lower bounds for the
  maximum number $\mu_{A_j}(d)$ of $A_j$-singularities on a surface of degree
  $d$ in most cases. 
  We also give a generalization to higher
  dimensions which leads to new lower bounds even in the case of nodal
  hypersurfaces in $\dP^n, n\ge5$.

  To conclude, we work out in detail a classical idea of
  B.~Segre which leads to some interesting examples, e.g.\ to a
  sextic with $36$ cusps.
\end{abstract}

\maketitle

%
%
%
\section{Introduction}
All possible configurations of singularities on a surface of degree $3$
in $\dP^3:=\dP^3(\dC)$ are known since Schläfli's work
\cite{Schl63} in the $19^\textup{th}$ century, see \cite{labsHolzerCS} and
\cite{labsAlgSurf} for explicit equations and illustrating pictures.  
In the case of degree $4$, the classification was completed recently by
Yang \cite{yangListQuartics} using computers. 

Much less is known for higher degrees, even when restricting to a particular
type of singularity.
E.g., the maximum number of $A_1$-singularities on a surface of degree $d$ is only
known for $d\le6$. 
We recently improved the case $d=7$ using computer algebra and geometry over
prime fields \cite{labs99}.   
The best lower bounds for surfaces of large degree $d$ with
$A_1$-singularities are given by Chmutov's construction \cite{chmuP3}.   

For higher singularities --- e.g., singularities of type $A_j$ which are locally
equivalent to $x^{j+1}+y^2+z^2$ --- the situation is even more difficult. 
We denote by $\mu_{A_j}(d)$ the maximum number of singularities of type $A_j$
a surface of degree $d$ in $\dP^3$ can have. 
Barth \cite{barQuint} constructed a quintic with $15$ singularities of type
$A_2$ (also called (ordinary) cusps), and the author constructed a sextic with
$35$ such singularities \cite{labsSext35} using computer algebra in
characteristic zero which showed $\mu_{A_2}(6)\ge35$. 
The detailed study of a generalization of an idea of B.~Segre \cite{BSegCon2}
which we give in the appendix leads to: $\mu_{A_2}(6)\ge36$ (see equations
(\ref{eqnF_1}) and (\ref{eqnApproxF_i}), and corollary
\ref{ThmSext36}). 
Recently, Barth and others considered the codes connected to surfaces with
three-divisible cusps in analogy to the codes related to even sets of nodes,
see \cite{barK39Cusps, barramsCC, ramsQuQuCusps}.   

In general, the best lower bounds for the $\mu_{A_j}(d)$ known up to now are
given by a direct generalization of an idea already used by Rohn in the
$19^\textup{\tiny th}$ century \cite[p.~33]{rohnNodQu}: $\mu_{A_j}(d) \ge
{1\over2}d(d-1)\lfloor {d\over j+1}\rfloor$ for
$d\ge 2(j+1)$ (see, e.g.,  \cite{barramsCC} and \cite{labsSext35} for
applications of this). 
For many degrees, one can also use the already mentioned generalization of
B.~Segre's construction \cite{BSegCon2} which is usually better than Rohn's if
it can be applied.

In the main part of the article we describe a variant of
Chmutov's construction \cite{chmuP3} which 
leads to the lower bound (corollary \vref{corMjd}):
\begin{equation}\mu_{A_j}(d) \ \gtrapprox \ {3j+2\over  
  6j(j+1)} d^3.\end{equation}
To our knowledge, this gives asymptotically the best known bounds for any
$j\ge2$. 
The construction reaches more than $\approx 75\%$ of the theoretical upper
bound in all cases. 
We compute this upper bound in section \vref{secUpperBounds}. 
Table \vref{tab_vch_largetable} gives an overview of our results for low $j$,
see also corollaries \ref{corMjd} and \ref{corQuotient}.
We describe a generalization of our construction to higher dimensions in
section \vref{secHighDim}. 
This leads to new lower bounds even in the case of nodal hypersurfaces. 

\begin{table}[htbp]
  \begin{center}
    \begin{tabular}{c|@{\,}c@{}c@{}c@{}c@{}c@{\ \,}c@{\ \,}c@{\ \,}c@{\
          \,}c@{\ \,}c@{\ \ }l} 
      \begin{tabular}{@{}c@{}}\setlength{\unitlength}{4144sp}%
\begingroup\makeatletter\ifx\SetFigFont\undefined%
\gdef\SetFigFont#1#2#3#4#5{%
  \reset@font\fontsize{#1}{#2pt}%
  \fontfamily{#3}\fontseries{#4}\fontshape{#5}%
  \selectfont}%
\fi\endgroup%
\begin{picture}(294,294)(-11,557)
\thinlines
{\color[rgb]{0,0,0}\put(  1,839){\line( 1,-1){270}}
}%
\put( 89,626){\makebox(0,0)[rb]{\smash{{\SetFigFont{9}{10.8}{\rmdefault}{\mddefault}{\updefault}{\color[rgb]{0,0,0}j}%
}}}}
\put(162,709){\makebox(0,0)[lb]{\smash{{\SetFigFont{9}{10.8}{\rmdefault}{\mddefault}{\updefault}{\color[rgb]{0,0,0}d}%
}}}}
\end{picture}\end{tabular} & $3$ & $4$ &
      $5$ & $6$ & $7$ & $8$ & $9$ & $10$ & $11$ & $12$ & $d$ \\ 
      \hline
      $1$ & \ublb{4}{4} & \ublb{16}{16} & \ublb{31}{31} & 
      \ublb{65}{65} & \ublb{99}{104} & \ublb{168}{174} & \ublb{216}{246} & 
      \ublb{345}{360} & \ublb{425}{480} & \ublb{600}{645} & $\approx \ \,
      \ublb{5/12}{4/9} \cdot d^3$\\

      $2$ & \ublb{3}{3} & \ublb{8}{8} & \ublb{15}{20} & \ublb{36}{37} & 
      \ublb{{52}}{62} & \ublb{{70}}{98} & \ublb{{126}}{144} & 
      \ublb{{159}}{202} & \ublb{{225}}{275} & \ublb{{300}}{363} &
      $\approx \, \ublb{{2/9}}{1/4} \cdot d^3$\\

      $3$ & \ublb{1}{1} & \ublb{6}{6} & \ublb{10}{13} & \ublb{15}{26} & 
      \ublb{{31}}{44} & \ublb{{64}}{69} & \ublb{{72}}{102} &
      \ublb{{114}}{144} & \ublb{{140}}{195} &
      \ublb{{198}}{258} & $\approx \ \ \; \ublb{{11/72}}{8/45} \ \cdot d^3$ \\

      $4$ & \ublb{1}{1} & \ublb{4}{4} & \ublb{10}{11} & \ublb{15}{20} & 
      \ublb{21}{35} & \ublb{{32}}{54} & \ublb{{54}}{80} & 
      \ublb{{100}}{112} & \ublb{{110}}{152} & 
      \ublb{{132}}{201} & $\approx \ \, \ublb{{7/60}}{5/36} \ \cdot d^3$ \\
    \end{tabular}
  \end{center}
  \caption{Known upper and lower bounds for the maximum
    number $\mu_{A_j}(d)$ of singularities of type $A_j, \ j=1,2,3,4,$ on a
    surface of degree $d$ in $\dP^3$. 
    For $j\ge2$ and $d\ge5$, the lower bounds are attained by our examples or
    by the generalization of B.~Segre's idea which we work out in the appendix.
  } 
  \label{tab_vch_largetable}
\end{table}

I thank D.\ van Straten for all his motivation, many valuable discussions, and for 
introducing me to the theory of Dessins d'Enfants.

%
%
%
\section{Chmutov's Idea} 

We start with some notation: A point $z_0\in\dC$ is a \emph{critical point of
  multiplicity} $j\in\dN$ 
of a polynomial $g\in\dC[z]$ in one variable if the first $j$ derivatives of
$g$ vanish at $z_0$: $g^{(1)}(z_0) = \cdots = g^{(j)}(z_0) = 0$.
The number $g(z_0)$ is called the \emph{critical value} of $z_0$.
A critical point of multiplicity $j, j>1,$ is called a \emph{degenerate}
critical point. 

In \cite{chmuP3}, Chmutov uses the following idea:
\begin{itemize}
\item
  Let $P_d(x,y) \in \dC[x,y]$ be a polynomial of degree $d$ with few different
  critical values, all of which are non-degenerate. 
  By a coordinate change, we may assume that the two critical values which
  occur most often are $0$ and $-1$. 
  We assume that they occur $\nu(0)$ and $\nu(-1)$ times, and that
  $\nu(0)>\nu(-1)$.  
\item
  Let $T_d(z) \in \dR[z]$ be the Chebychev polynomial of degree $d$ with
  critical values $-1$ and $+1$, where $-1$ occurs $\lfloor {d\over2}\rfloor$
  times and $+1$ occurs $\lfloor {d-1\over2}\rfloor$ times.
\item
  It is easy to see that the projective surface given by the affine equation
  \begin{equation}\label{chmuOrig}
    P_d(x,y) + \frac{1}{2}(T_d(z)+1) = 0
  \end{equation}
  has $\nu(0) \cdot
  \lfloor {d\over2}\rfloor + \nu(-1) \cdot \lfloor {d-1\over2}\rfloor$ nodes.
\end{itemize}

Chmutov uses for $P_d(x,y)$ the so-called folding polynomials $F^{A_2}_d$
associated to the root system $A_2$ (see \cite{witFoldPoly}): 
{\small 
  \begin{equation}\label{eqnFA2d}
    F^{A_2}_d(x,y) := 2+\det\left(
      \begin{array}{c@{\,}c@{\,}c@{\,}c@{\,}c@{\,}c@{\,}c}
        \rule{0pt}{1.4em}x & 1 & 0 & \cdots & \cdots & \cdots & 0\\[-0.6em]
        \rule{0pt}{1.4em}2y & x & \ddots & \ddots &  &  & \vdots \\[-0.6em]
        \rule{0pt}{1.4em}3 & y & \ddots & \ddots & \ddots & & \vdots \\[-0.6em]
        \rule{0pt}{1.4em}0 & 1 & \ddots & \ddots & \ddots & \ddots & \vdots \\[-0.6em]
        \rule{0pt}{1.4em}\vdots & \ddots & \ddots & \ddots & \ddots & \ddots & 0\\[-0.6em]
        \rule{0pt}{1.4em}\vdots &  & \ddots & \ddots & \ddots & \ddots & 1\\[-0.6em]
        \rule{0pt}{1.4em}0 & \cdots & \cdots & 0 & 1 & y & x
      \end{array}
    \right) + 
    \det\left(
      \begin{array}{c@{\,}c@{\,}c@{\,}c@{\,}c@{\,}c@{\,}c}
        \rule{0pt}{1.4em}y & 1 & 0 & \cdots & \cdots & \cdots & 0\\[-0.6em]
        \rule{0pt}{1.4em}2x & y & \ddots & \ddots &  &  & \vdots \\[-0.6em]
        \rule{0pt}{1.4em}3 & x & \ddots & \ddots & \ddots & & \vdots \\[-0.6em]
        \rule{0pt}{1.4em}0 & 1 & \ddots & \ddots & \ddots & \ddots & \vdots \\[-0.6em]
        \rule{0pt}{1.4em}\vdots & \ddots & \ddots & \ddots & \ddots & \ddots & 0\\[-0.6em]
        \rule{0pt}{1.4em}\vdots &  & \ddots & \ddots & \ddots & \ddots & 1\\[-0.6em]
        \rule{0pt}{1.4em}0 & \cdots & \cdots & 0 & 1 & x & y
      \end{array}
    \right).\end{equation}
}

\begin{figure}[htbp] 
\begin{center}
  \begin{tabular}{ccc}
    \includegraphics[width=1in]{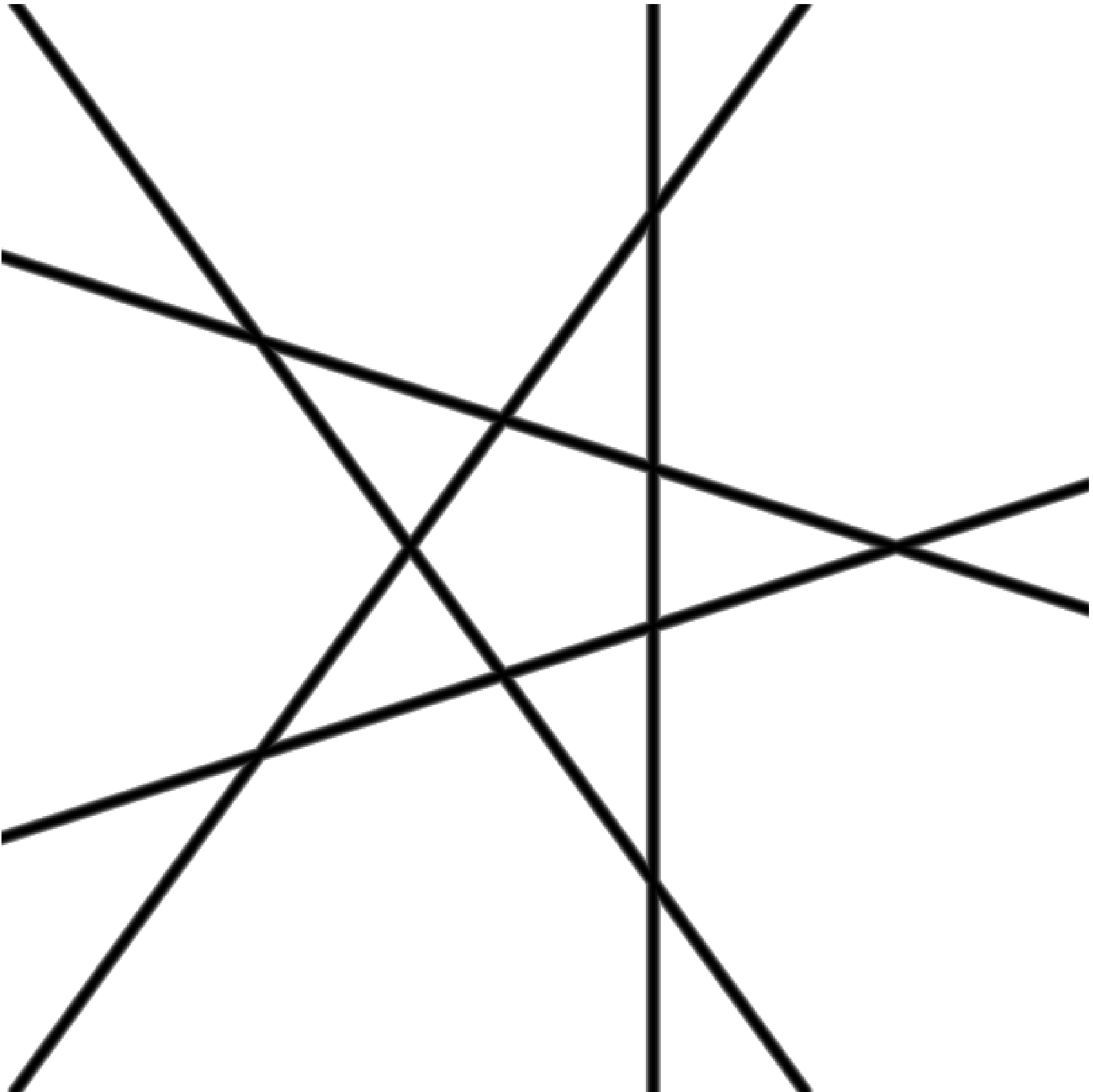} & 
    \includegraphics[width=1in]{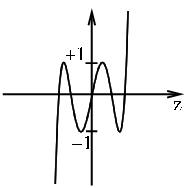} & 
    \includegraphics[width=1in]{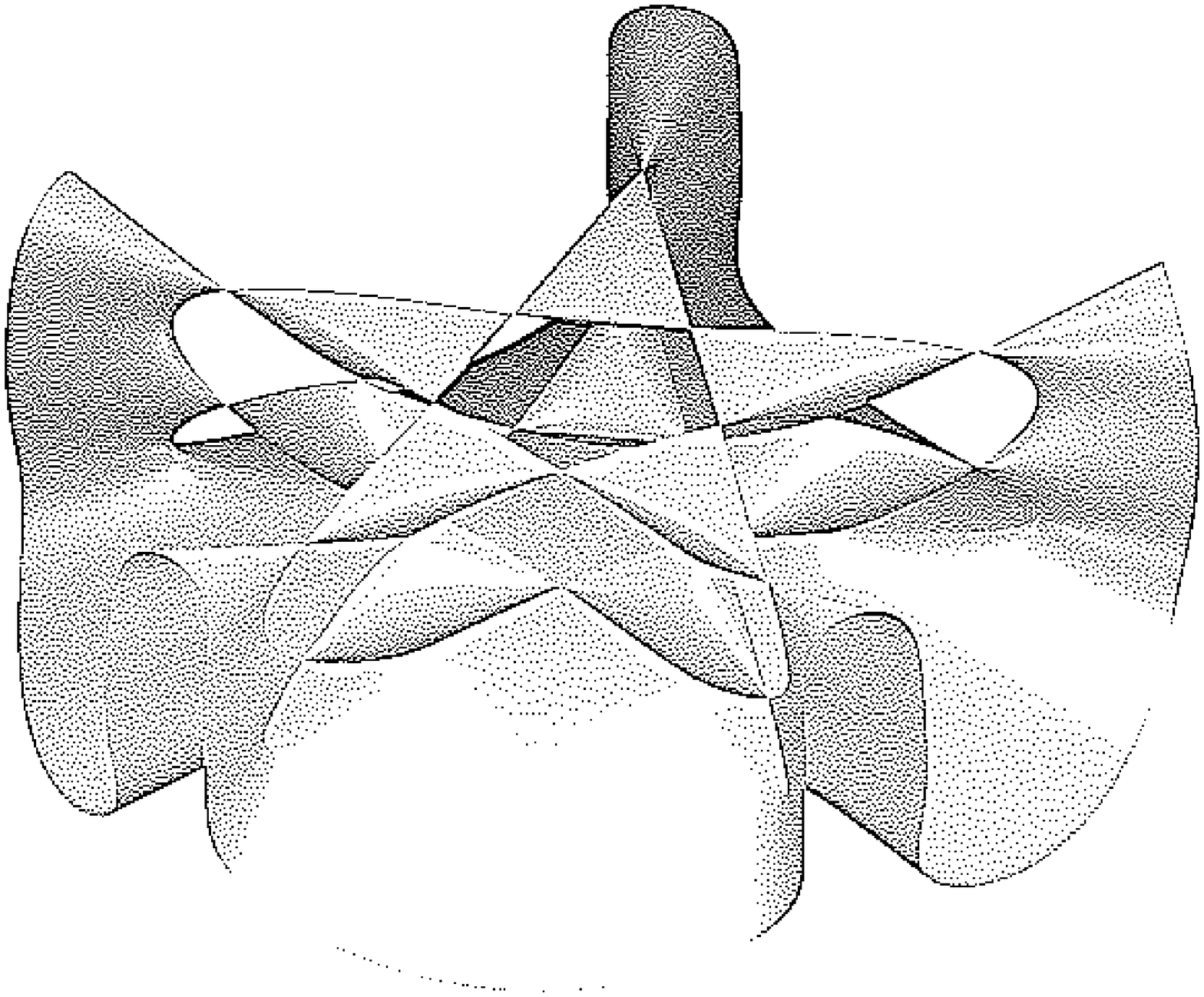}\\
    $R_5(x,y)$ & 
    $T_5(z)$ & 
    $R_5(x,y) + \frac{1}{2}(T_5(z)+1)$
  \end{tabular}
  \caption{A variant of Givental's and Chmutov's construction: A regular
    $5$-gon 
    $R_5(x,y)$, the Chebychev polynomial $T_5(z)$ 
    and the surface $R_5(x,y) + \frac{1}{2}(T_5(z)+1)$ with $10\cdot2 +
    5\cdot 2 = 30$ nodes.} 
  \label{figDucoQuint}
\end{center}
\end{figure}

$F^{A_2}_d(x,y)$ has $d\choose2$ critical points with critical
value $0$ and ${1\over3}d(d-3)$ 
critical points with critical value $-1$ if $d \equiv 0 \mod
3$, and ${1\over3}(d(d-3)+2)$ otherwise (see \cite{chmuP3}); 
the other critical points have critical value $8$.
To our knowledge, these are still the best known polynomials for this
purpose; in \cite{chmuCritVals}, Chmutov conjectured them to be asymptotically
the best.
We illustrate the idea using a variant of Chmutov's construction which was
suggested in the case of cubic hypersurfaces by Givental \cite[p.\ 419]{AGV2}: 
We take a regular $d$-gon $R_d(x,y)$ for $P_d(x,y)$ (see fig.\
\vref{figDucoQuint}).  
This has $d\choose2$ critical points with critical value $0$, one critical
point over the origin, and $d$ critical points with each of the other critical
values (we assume that one of these is $-1$).  
Then the construction above gives $30$ nodes for $d=5$, see fig.\
\vref{figDucoQuint}. 
Notice that $F^{A_2}_5$ only leads to $28$ nodes, but for all $d\ge6$,
$F^{A_2}_d$ is much better than $R_d$.

%
%
%
\section{Adaption to Higher Singularities}
\label{secAdapt}

To adapt Chmutov's construction (\ref{chmuOrig}) to higher singularities of
type $A_j$, we replace the polynomials $T_d(z)$ by polynomials with degenerate
critical points.  

For the construction of a quintic surface with many cusps, we thus take again
the regular $5$-gon $R_5(x,y)\in\dR[x,y]$ together with a polynomial
$T^2_5(z)\in\dR[z]$ of degree $5$ with the maximum number of critical points
of multiplicity two.  
\begin{figure}[htbp]
\begin{center}
  \begin{tabular}{ccc}
    \includegraphics[width=1in]{n_gon_5} & 
    \includegraphics[width=1in]{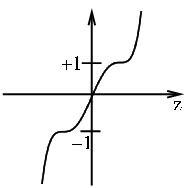} & 
    \includegraphics[width=1in]{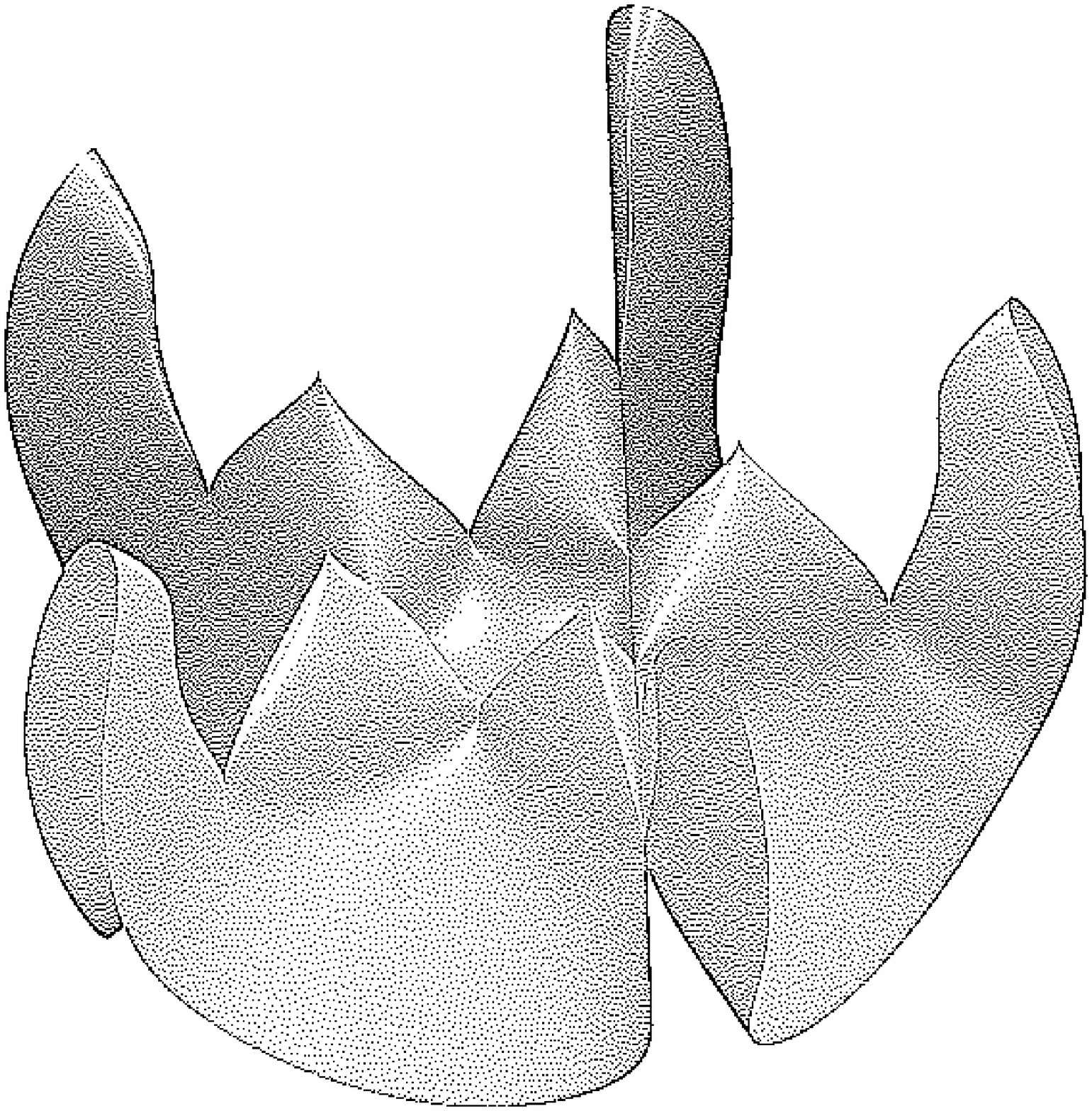}\\[0.3em]
    $R_5(x,y)$ & 
    $T^2_5(z)$ & 
    $R_5(x,y) + \frac{1}{2}(T^2_5(z)+1)$ 
  \end{tabular}
  \caption{The construction of a quintic with $15$ cusps.} 
  \label{figA2Quint}
\end{center}
\end{figure}  
As the derivative of such a polynomial has degree $4$, the maximum number of
such critical points is ${4\over2}=2$, see fig.\ \vref{figA2Quint}.
The critical values of these two critical points have to be different because
a horizontal line through both critical points would intersect the curve in
six points counted with multiplicities.
Similar to the situation for nodes in (\ref{chmuOrig}) the surface
$R_5(x,y) + \frac{1}{2}(T^2_5(z)+1)$ has $10\smcdot1 + 5\smcdot1=15$
singularities of type $A_2$.

As mentioned in the introduction, Barth already constructed another quintic
with $15$ cusps \cite{barQuint}. 
The author constructed a sextic with $35$ cusps in \cite{labsSext35}, and the
appendix gives a sextic with $36$ cusps.
But our variant of Chmutov's construction which will be presented in the
following sections gives new lower
bounds for the maximum number of $A_2$-singularities for all degrees $d\ge7$. 
We take surfaces in separated variables defined by polynomials of the form: 
\begin{equation}\label{eqnChjd}\textup{Chm}(G^j_d) := F^{A_2}_d + G^j_d,\end{equation} 
where $F^{A_2}_d(x,y)\in\dR[x,y]$ is the folding polynomial defined in
(\ref{eqnFA2d}) and where $G^j_d(z)\in\dC[z]$ is a polynomial of 
degree $d$ with many critical points of multiplicity $j$ with critical values
$-1$ and $+1$. 
E.g., for $j=1$, the ordinary Chebychev polynomials $G^1_d(z) := T_d(z)$ yield
to Chmutov's surfaces with many nodes.
In the following sections, we discuss two generalizations of the ordinary
Chebychev polynomials to polynomials with critical points of higher
multiplicity which give surfaces of degree $d$ with many $A_j$-singularities,
$j < d$.

%
%
%
\section{$j$-Belyi Polynomials via Dessins d'Enfants}
\label{secDdE}

The existence of polynomials in one variable with only two different critical
values with prescribed multiplicities of the critical points can be
established using ideas of Hurwitz \cite{hurVerzwPkt} based on Riemann's
Existence Theorem.    
The interest in this subject was renewed by Grothendieck's \emph{Esquisse d'un
  programme}. 
Nowadays, it is commonly known under the name of \emph{Dessins d'Enfants}.  
We will use the following proposition / definition which is basically taken from
\cite{adrzvoCompPlaneTrees}:
\begin{PropDef} \ 
  \begin{enumerate}
  \item
    A tree (i.e.\ a graph without cycles) with a prescribed cyclic order of the edges
    adjacent to each vertex is called a {\bf plane tree}.
    A plane tree has a natural bicoloring of the vertices (black/white). 
    If we fix the color of one vertex, then this bicoloring is unique.
  \item
    A polynomial with not more than two different critical values is
    called a {\bf Belyi polynomial}. 
  \item
    For a given Belyi polynomial $p:\dC\to\dC$ with critical values $c_1$ and $c_2$, we
    define the {\bf plane tree $PT(p)$ associated to $p$} to be the inverse
    image $p^{-1}([c_1,c_2])$ of the 
    interval $[c_1,c_2]$, where $p^{-1}(c_1)$ are the black vertices, and
    $p^{-1}(c_2)$ are the white vertices of the tree (see fig.\ \vref{figDessinsProj}).
    \begin{figure}[htbp]
      \begin{center}
        \begin{tabular}{c@{\qquad}c}
          \begin{tabular}{c}\includegraphics[width=1.5in]{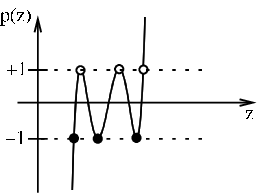}\end{tabular} & 
          \begin{tabular}{c}\includegraphics[width=1.5in]{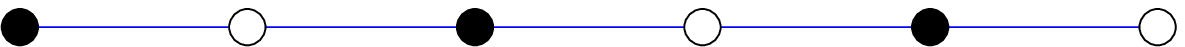}\end{tabular} 
        \end{tabular}
        \caption{The ordinary Chebychev polynomial $T_5$ with two critical points
          with critical value $-1$ and two with critical value $+1$. 
          The right picture shows its plane tree $PT(T_5)$.
          A vertex with two adjacent edges corresponds to a critical point with
          multiplicity $1$, a vertex with one adjacent edge corresponds to a
          non-critical point.} 
        \label{figDessinsProj}
      \end{center}
    \end{figure} 
  \item
    For any plane tree, there exists a Belyi polynomial whose
    critical points have the multiplicities given by the number of edges
    adjacent to the vertices minus one and vice verca. 
  \end{enumerate}
\end{PropDef}

We will need the following two trivial bounds concerning critical points:
\begin{Lem}\label{lem_gub}
  Let $d,j\in\dN$.
  Let $g\in\dC[z]$ be a polynomial of degree $d$ in one variable with only isolated
  critical points. Then:
  \begin{enumerate}
    \item
      The total number of different critical points of $g$ of multiplicity $j$
      does not exceed $\lfloor {d-1\over j} \rfloor$. 
    \item 
      The number of different critical points of $g$ of multiplicity $j$ with
      the same critical value does not exceed $\lfloor {d\over j+1}
      \rfloor$.\hfill$\Box$ 
  \end{enumerate}
\end{Lem}

We give a special name to polynomials reaching the first of these bounds:
\begin{Def}
  Let $d,j\in\dN$ and let $p$ be a Belyi polynomial of degree $d$. 
  We call $p$ a {\bf $j$-Belyi polynomial} if $p$ has the maximum possible
  number $\lfloor {d-1\over j}\rfloor$ of critical points of multiplicity
  $j$.
\end{Def}

\begin{Exa}
  The ordinary Chebychev polynomials $T^1_d(z) := T_d(z)$ are $1$-Belyi polynomials. 
  $T^2_5(z)$ in fig.\ \vref{figA2Quint} is a $2$-Belyi
  Polynomial.\hfill$\Box$
\end{Exa}

A special type of $j$-Belyi polynomials are those of degree $j+1$. 
We will join several plane trees corresponding to such $j$-Belyi polynomials
of degree $j+1$ to form larger plane trees in the following sections:
\begin{Def}
  We call the plane tree corresponding to a $j$-Belyi polynomial of degree
  $j+1$ a {\bf $j$-star}. 
  If the center of this tree is a black (resp.\ white) vertex we call it a
  $\bullet$- (resp.\ $\circ$-) {\bf centered $j$-star} (see fig.\
  \vref{figJStar}). 
  \begin{figure}[htbp]
    \begin{center}
      \begin{tabular}{c@{\qquad\qquad}c}
        \begin{tabular}{c}\includegraphics[width=1in]{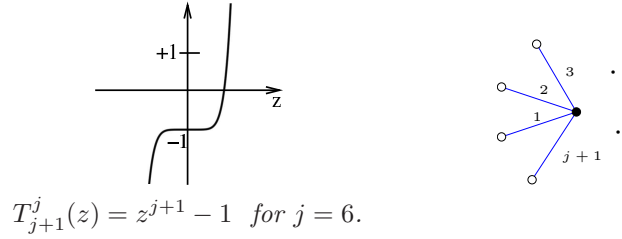}\\
          $T^j_{j+1}(z) = z^{j+1}-1$ \ for $j=6$.\end{tabular} & 
      \begin{tabular}{c}\input{St_n.pstex_t}\end{tabular} 
      \end{tabular} 
      \caption{The polynomial
        $T^j_{j+1}(z)$ with exactly one critical point $z_0=0$ of multiplicity
        $j$ and critical value $-1$ together with the corresponding
        $\bullet$-centered $j$-star.}   
      \label{figJStar}
    \end{center}  
  \end{figure}
\end{Def}

%
%
%
\section{The Polynomials $T^j_d(z)$}
\label{secTjd}
\label{RET}

A natural generalization of the ordinary Chebychev polynomials to polynomials
$G^j_d(z)$ with degenerate critical points that can be used in the
construction of equation (\ref{eqnChjd}) on page \pageref{eqnChjd} comes from
the following intuitive idea: 
Take polynomials which look similar to the ordinary Chebychev polynomials (fig.\
\ref{figDessinsProj}), but which have higher vanishing derivatives such that
they are $j$-Belyi polynomials. 
\begin{Exa}
  A $3$-Belyi polynomial of degree $13$ has
  $\bigl\lfloor\frac{13-1}{3}\bigr\rfloor=4$ critical points of 
  multiplicity $3$. 
  The polynomial $T^3_{13}$ has two critical points with critical value $-1$
  and two with critical value $+1$.
  The plane tree showing the existence of such a polynomial consists of four 
  connected $3$-stars. 
  To show the similarity to the ordinary Chebychev polynomials we draw them in
  fig.\ \ref{figTjdTT} as four bouquets of $1$-stars attached to the plane
  tree in fig.\ \vref{figDessinsProj}.
  A straightforward {\sc Singular} \cite{Singular} script to
  compute the equation of $T^3_{13}(z)$ can be found on the website
  \cite{labsAlgSurf}.  
  \hfill$\Box$
  \begin{figure}[htbp]
    \begin{center}
      \begin{tabular}{cc}
        \\[0.0em]
      \begin{tabular}{c}
        \includegraphics[height=0.8in]{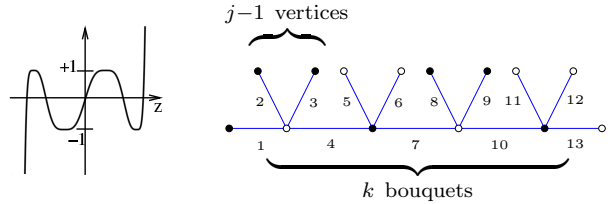}
      \end{tabular} &
      \begin{tabular}{c}
        \input{plane_trees_shabat_A3_numbered.pstex_t} 
      \end{tabular} 
      \end{tabular}\\[0.2em]
      \caption{The bicolored plane tree $PT(T^j_d)$ for the polynomial $T^j_d(z)$ for 
        $j=3$, $d=13$, $k:={d-1\over j}=4$. 
        It consists of $k$ connected $j$-stars.
        Here, we line them up to show the similarity to the ordinary Chebychev
        polynomials in fig.\ \vref{figDessinsProj}. 
        See \cite{labsAlgSurf} for a {\sc Singular} \cite{Singular} script to
        compute the equation of $T^3_{13}(z)$.        
      }  
      \label{figTjdTT}
    \end{center}
  \end{figure} 
\end{Exa}

\begin{ThmDef}
  Let $d,j\in\dN$ with $d>j$.
  There exists a polynomial $\mathbf{T^j_d(z)}$ of degree $d$ with
  $\lceil{1\over2}\lfloor{d-1\over j}\rfloor\rceil$ critical points of
  multiplicity $j$ with critical value $-1$ and
  $\lfloor{1\over2}\lfloor{d-1\over j}\rfloor\rfloor$ such critical points
  with critical value $+1$. 
\end{ThmDef}
\begin{proof}
  The corresponding plane tree $PT(T^j_d)$ can be defined as follows (compare
  fig.\ \vref{figTjdTT}). 
  For $d=k\cdot (j+1), \ k\in\dN$, we take $k$ connected $j$-stars.
  Fixing the center of the first $j$-star to be white, the plane tree has a 
  unique bicoloring.   
  If $d=l+k\cdot (j+1)$ for some $1\le l\le j$, we attach another $l$-star
  to get a polynomial of degree $d$.
\end{proof}

Although there is an explicit recursive construction of ordinary Chebychev
polynomials and their generalizations to higher dimensions (so-called folding
polynomials, see \cite{witFoldPoly}), we do not know a similar explicit
construction of the polynomials $T^j_d(z)$ for $j\ge2$.   
To our knowledge, they can only be computed for low degree $d$ until now,
e.g.\ using Groebner Basis.  
When plugged into the construction (\ref{eqnChjd}) on page \pageref{eqnChjd}
the existence of the polynomials $T^j_d$ immediately implies: 
\begin{Cor}\label{corTjd}
  Let $d,j\in\dN$ with $d>j$.  
  There exist surfaces 
  $$\textup{Chm}(T^j_d) := F^{A_2}_d + \frac{1}{2}(T^j_d+1)$$ 
  of degree $d$ with the following number of singularities of type $A_j$: 
  $$\qquad\begin{array}{ll}
    {1\over2}d(d-1) \smcdot \lceil{1\over2}\lfloor{d-1\over j}\rfloor\rceil 
    + {1\over3}d(d-3) \smcdot \lfloor{1\over2}\lfloor{d-1\over j}\rfloor\rfloor, 
    & \textup{if} \ d \equiv 0 \mod 3,\\[0.5em]
    {1\over2}d(d-1) \smcdot 
    \lceil{1\over2}\lfloor{d-1\over j}\rfloor\rceil  +
    {1\over3}(d(d-3)+2) \smcdot \lfloor{1\over2}\lfloor{d-1\over j}\rfloor\rfloor
    & \textup{otherwise}.\qquad\qquad\qquad\Box 
  \end{array}$$
\end{Cor}

%
%
%
\section{The Polynomials $M^j_d(z)$}
\label{secMjd}

The $j$-Belyi polynomials $T^j_d(z)$ described in the previous section reach
the first bound of lemma \vref{lem_gub}. 
The $j$-Belyi polynomials $M^j_d(z)$ whose existence will be shown in this section  
also achieve the second bound of this lemma. 
We start with two examples:

\begin{Exa}
  The $2$-Belyi polynomial $T^2_9(z)$ is the example of the smallest degree
  from the previous section that does not reach the second bound of lemma
  \ref{lem_gub}.  
  The plane tree $PT(M^2_9(z))$ in fig.\ \vref{figM29} shows the existence of a
  $2$-Belyi polynomial of degree $9$ that achieves this bound.

  As in the case of the polynomials $T^j_d(z)$, it is possible to compute the
  polynomials $M^j_d(z)$ explicitly for low $j$ and $d$.
  For our case $j=2, d=9$ we denote by $u$ the unique critical point with
  critical value $+1$ and by $b_0, b_1, b_2$ the three critical points with
  critical value $-1$. 
  When requiring $b_2=0$ (i.e., $M^2_9(0)=-1$), $M^2_9(z)$ has the derivative 
  $$\frac{\partial M^2_9}{\partial z}(z) = (z-b_0)^2\cdot (z-b_1)^2\cdot
  z^2\cdot (z-u)^2.$$ 
  Using {\sc Singular} \cite{Singular}, we find: $u^9=18$ and $b_0$
  and $b_1$ are the two distinct roots of $z^2-3uz+3u^2=0$.
  Notice that $b_0, b_1 \notin \dR$ even if we take $u\in\dR$.
  \hfill$\Box$ 
  \begin{figure}[htbp]
    \begin{center}
      \begin{tabular}{ccccc}
        \includegraphics[scale=0.28]{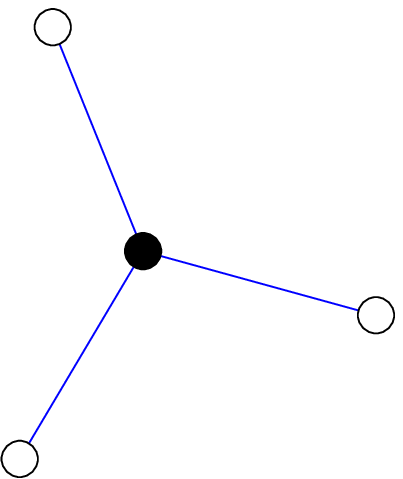} & 
        \quad\quad &
        \includegraphics[scale=0.28]{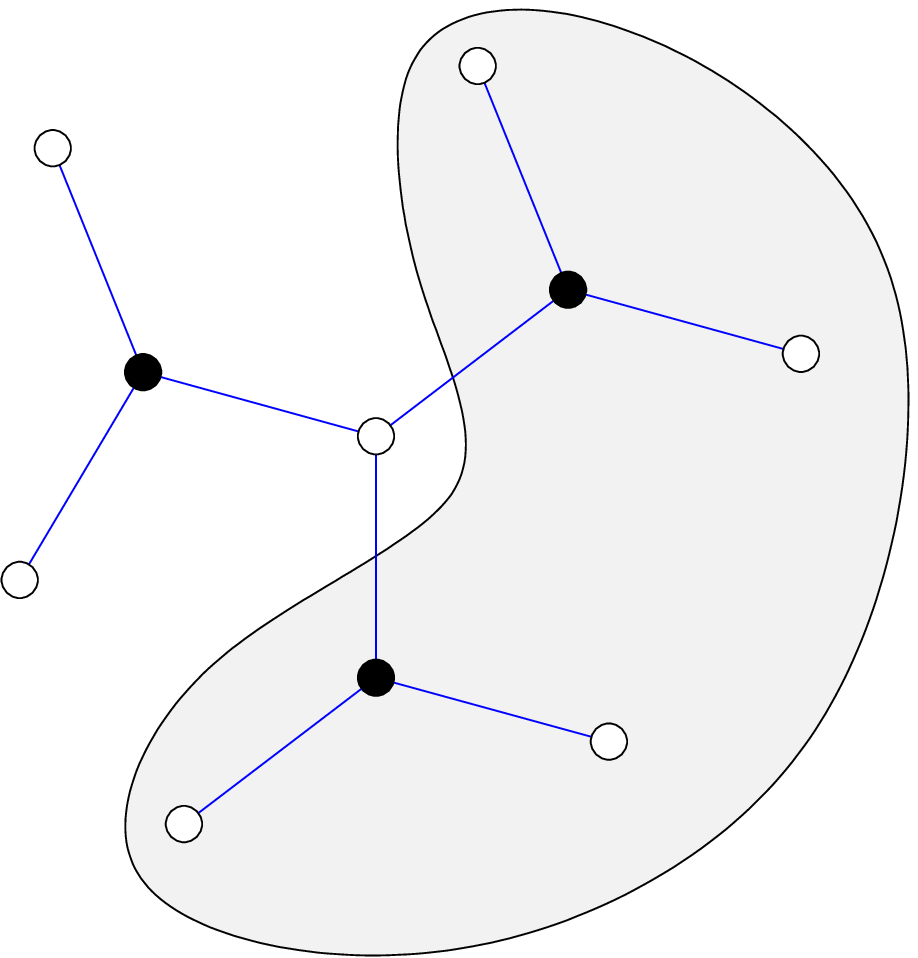} & 
        \quad\quad &
        \includegraphics[scale=0.28]{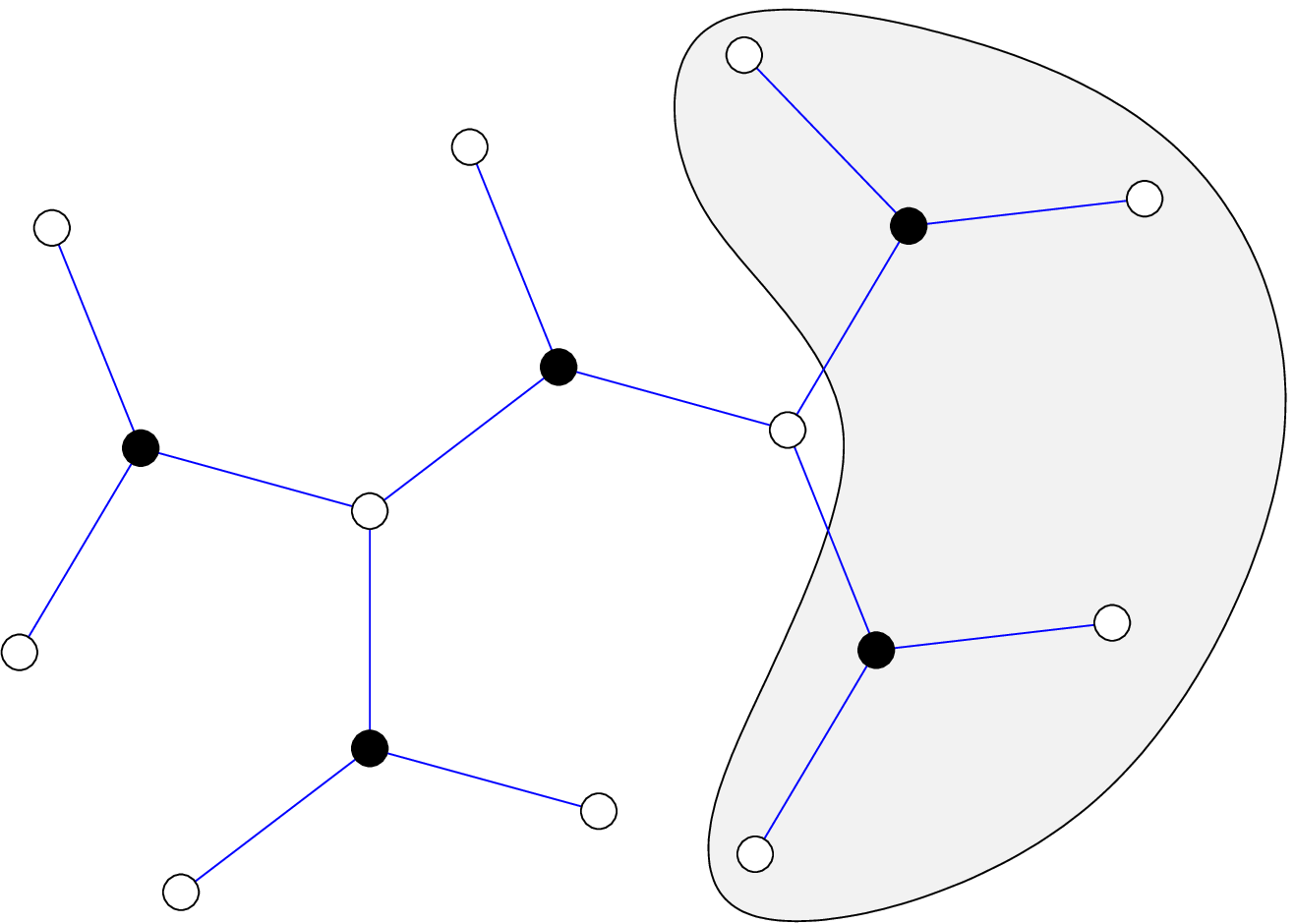}\\[0.5em]
        (a) \ $PT(M^2_3)$ & & 
        (b) \ $PT(M^2_9)$ & & 
        (c) \ $PT(M^2_{15})$ 
      \end{tabular}
      \caption{
        To obtain $PT(M^2_9)$ from the $2$-star $PT(M^2_3)=PT(T^2_3)$, we attach two
        $\bullet$-centered $2$-stars to one of the $\circ$-vertices (marked by
        the grey background). 
        The corresponding polynomial $M^2_9(z)$ has thus $3$ critical points of multiplicity
        $2$ with critical value $-1$ (the $3$ $\bullet$-centered $2$-stars) and
        $1$ such point with critical value $+1$ (the only $\circ$-centered
        $2$-star). 
        $M^2_{15}$ has five and two, respectively.
      } 
      \label{figM29}
    \end{center}
  \end{figure} 
\end{Exa}

\begin{Exa}
  If $d\ne k \cdot (j+1)$ for some $k\in\dN$, the
  construction of a plane tree corresponding to a polynomial reaching both
  bounds of lemma \ref{lem_gub} is a little more delicate than in the previous
  example.
  The cases $PT(M^2_{11})$ and $PT(M^2_{12})$ in fig.\ \vref{figM3}
  illustrate this.\hfill$\Box$
  \begin{figure}[htbp]
    \begin{center}
      \begin{tabular}{ccc}
        \includegraphics[scale=0.28]{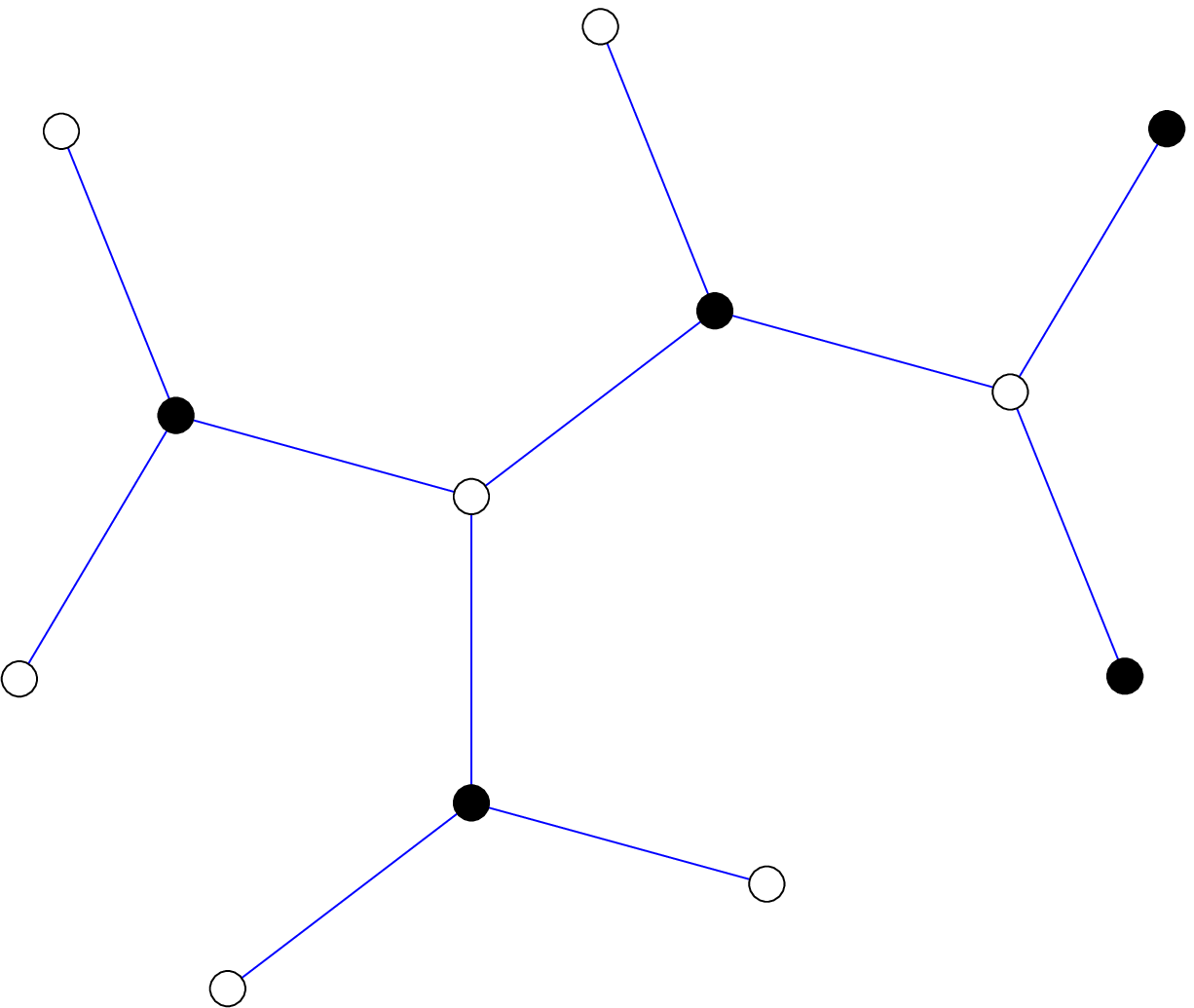} & 
        \quad\quad &
        \includegraphics[scale=0.28]{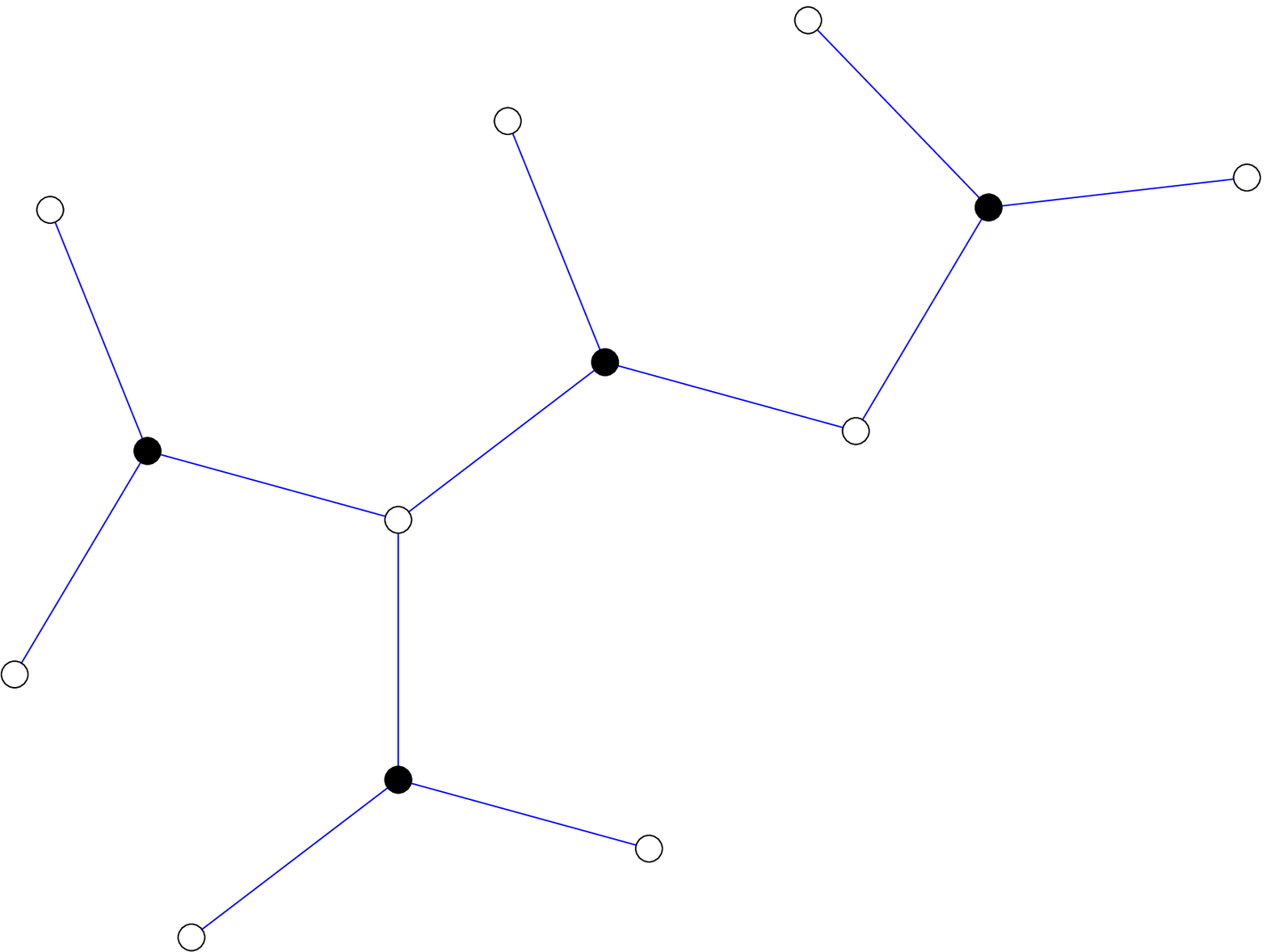}\\[0.5em]
        (a) \ $PT(M^2_{11})$ & & 
        (b) \ $PT(M^2_{12})$ 
      \end{tabular}
      \caption{
        $M^2_{11}$ and $M^2_{12}$ have the same number of critical points of
        multiplicity $j$.
        $M^2_{12}$ has five ones with critical value $-1$ and only one
        critical point with critical value $+1$. 
        $M^2_{11}$ has three critical points with critical value $-1$
        and two with critical value $+1$. 
      } 
      \label{figM3}
    \end{center}
  \end{figure} 
\end{Exa}

\begin{ThmDef}
  Let $d,j\in\dN$ with $d>j$.
  There exists a polynomial $\mathbf{M^j_d(z)}$ of degree $d$ with
  $\left\lfloor{d\over j+1}\right\rfloor$ critical points of
  multiplicity $j$ with critical value $-1$ and
  $\left(\left\lfloor{d-1\over
      j}\right\rfloor-\left\lfloor{d\over 
      j+1}\right\rfloor\right)$ such critical points
  with critical value $+1$. 
\end{ThmDef}
\begin{proof}
  The existence of a corresponding plane tree $PT(M^j_d)$ can be shown as
  follows (compare fig.\ \vref{figM29}). 
  For $d=j+1$ we define $PT(M^j_d)$ as a $\bullet$-centered $j$-star.
  For $d=(j+1) + k\cdot j\smcdot (j+1), \ k\in\dN$, we attach successively 
  sets of $j$ $\bullet$-centered $j$-stars as illustrated in figure
  \ref{figM29}.    
  If $d\ne (j+1) + k\cdot j\smcdot (j+1)$ for some $k\in\dN$ the existence of
  plane trees $PT(M^j_d)$ can be shown similarly (see fig.\ \vref{figM3}).
\end{proof}

The existence of the polynomials $M^j_d(z)$ has two immediate consequences:
\begin{Cor}\label{corBoundSharp}
  The bounds in lemma \vref{lem_gub} are sharp.\hfill$\Box$
\end{Cor}

It is clear that the polynomials $M^j_d$ cannot have only real coefficients
and only real critical points for $d$ large enough.
So, the same holds for the singularities of the surfaces of the following
corollary: 

\begin{Cor}\label{corMjd}
  Let $d,j\in\dN$ with $d>j$.  
  There exist surfaces $$\textup{Chm}(M^j_d) := F^{A_2}_d + M^j_d$$ of
  degree $d$ with the following number of singularities of type $A_j$: 
  $$\begin{array}{ll}
    {1\over2}d(d-1) \smcdot \left\lfloor{d\over j+1}\right\rfloor
    + {1\over3}d(d-3) \smcdot \left(\left\lfloor{d-1\over
      j}\right\rfloor-\left\lfloor{d\over 
      j+1}\right\rfloor\right), 
    & \textup{if} \ d \equiv 0 \mod 3,\\[0.5em]
    {1\over2}d(d-1) \smcdot \left\lfloor{d\over j+1}\right\rfloor +
    {1\over3}(d(d-3)+2) \smcdot \left(\left\lfloor{d-1\over
      j}\right\rfloor-\left\lfloor{d\over 
      j+1}\right\rfloor\right)
    & \textup{otherwise}.\qquad\qquad\qquad\Box 
  \end{array}$$
\end{Cor}

%
%
%
\section{Upper Bounds}
\label{secUpperBounds}

To get an idea of the quality of our best lower bounds given by our examples
$\textup{Chm}(M^j_d)$ from corollary \vref{corMjd} we compare them with
the best known upper bounds: Miyaoka's bound \cite{miyp3} and Varchenko's
Spectral Bound \cite{varBound}. 

\subsection{Varchenko's Spectral Bound}

It is well-known that Varchenko's Spectral Bound 
$\textup{Var}_{A_j}(d)$ for the maximum number of $A_j$-singularities \cite{varBound}
is not as good as Miyaoka's bound for fixed $j$ and large $d$. 
Although it is known that both bounds can be described by a polynomial of
degree three in $d$, we could not find explicit statements for Varchenko's
bound for $j>1$  in the literature. 
So, we compute these polynomials here using a short {\sc Singular}
\cite{Singular} script.   
The code can be downloaded from \cite{labsAlgSurf}.
In the following, we explain briefly how the algorithm works.

For even degree $d\ge4$ the spectrum $sp(d)$ of the singularity $x^d + y^d + 
z^d=0$ in $\dC^3$ consists of the spectral numbers $s_d(i) = {i+2\over d},
i=1,2,\dots,3(d-1)-2$, with multiplicities $m_d(i)$, 
where 
\begin{itemize}
\item
$m_d(1)=1$,
\item
$m_d(i+1) = m_d(i)+1+i, \quad i<d-1$,
\item
$m_d(i+1) = m_d(i)+2(i_{mid}-i) + 1, \quad d-1 \le i < i_{mid} := {3d\over 2}-2$, 
\item 
$m_d(3(d-1)-1-i) = m_d(i), \quad 1 \le i \le i_{mid}$ (symmetry of the spectrum).
\end{itemize} 
The spectrum of an $A_j$ singularity is also well-known (see e.g.\
\cite[p.\ 389]{AGV2}). 
Its spectral numbers are ${j+2 \over j+1}$, ${j+3 \over j+1}$, \dots, ${2j+1
  \over j+1}$, all with multiplicity $1$.
\begin{Exa}\label{exa62}
  The spectrum $sp(6)$ of the singularity $x^6+y^6+z^6$ is:
  \begin{center}
    \begin{tabular}{c}
      $\begin{array}{rc@{\ }c@{\ }c@{\ }c@{\ }c@{\ }c@{\ }c@{\ }c@{\ }c@{\ }c@{\ }c@{\
          }c@{\ }c@{\ }c@{\ }c@{\ }c@{\ }c@{\ }c@{\ }c@{\ }c@{\ }c@{\ }c@{\ }c@{\
          }c@{\ }c@{\ }c@{\ }c@{\ }c@{\ }c@{\ }c@{\ }c@{\ }c@{\ }c@{\ }c@{\ }c}  
        \footnotesize i& 1 && 2 && 3 && 4 && 5 && 6 && 7 && 8 && 9 && 10 && 11
          && 12 && 13\\[0.2em]
        \footnotesize \textup{spectral number} \ s_i& 
        3 \over 6& &4 \over 6& &5 \over 6& &6 \over 6& &7 \over 6& &8 \over 6& &9
        \over 6& &10 \over 6& &11 \over 6& &12 \over 6& &13 \over 6& &14 \over 6& &15
        \over 6\\[0.2cm]   
        \footnotesize \textup{multiplicity} \ m_i& 
        1& &3& &6& & 10& &15& &18& &19& &18& &15& &10& &6& &3& &1\\[0.1cm] 
      \end{array}$
    \end{tabular}
  \end{center}
  The spectral numbers of the $A_2$-singularity are: $\frac{8}{6}, \frac{10}{6}$,
  both with multiplicity $1$.\hfill$\Box$  
\end{Exa}

To compute Varchenko's bound we have to choose an open interval
$I = (\frac{i_l+2}{d}, \frac{i_r+2}{d})$ of length $1$ 
of the spectrum $sp(d)$ that contains all spectral numbers of the $A_j$
singularity and such that the sum of the multiplicities of the spectral
numbers in the interval is minimal.
Then we have to sum up all the multiplicities in this interval and divide by
$j$.

Let us write $d = k\smcdot(j+1) + l$. 
Then we may choose $I := (\frac{i_r+2-d}{d}, \frac{i_r+2}{d})$, where $i_r :=
k\smcdot(2j+1) + \left\lfloor{l\cdot(2j+1) \over j+1}\right\rfloor - 1$. 
We introduce some notations: $n_l := i_{mid}-(d-1)$, $n_r := i_r - i_{mid} - 1$, $n_{ll}
:= d-1 - n_l - n_r$, $m_{mid} = \sum_{i=1}^{d-1}i + (\frac{d}{2}-1)^2$. 
Using these we can compute Varchenko's bound $\textup{Var}_{A_j}(d)$ for the maximum
number of $A_j$-singularities on a surface of degree $d$ in $\dP^3$ for the
case $d,j\in\dN$ with $d\ge4$: 
  \begin{equation}\label{eqnVarB}
    \quad\begin{array}{r@{\ }c@{\ }l}
      \textup{Var}_{A_j}(d) &=& \Bigl\lfloor\frac{1}{j}\smcdot\Bigl(
        \frac{1}{2}\smcdot\bigl(
          \sum_{i=1}^{d-1}i \ + \sum_{i=1}^{d-1}i^2 \ - \sum_{i=1}^{d-1-n_{ll}}i \
          - \sum_{i=1}^{d-1-n_{ll}}i^2 \bigr)\\
      &&    + (n_r + n_{ll})\cdot m_{mid} \ - \sum_{i=1}^{n_r} i^2 \ -
        \sum_{i=1}^{n_l-1} i^2 
      \Bigr)\Bigr\rfloor.\quad\qquad\qquad\Box
    \end{array}
  \end{equation}

\begin{Exa}
  Let us look at the case $d=6, j=2$ as in example \vref{exa62}.
  In this case, the constants used above have the following values: 
  $k=2$, $l=0$, $i_r=9$, $i_l=3$, $i_{mid}=7$, $n_l = 2$, $n_r = 1$, $n_{ll} =
  2$, $m_{mid} = 19$. 
  We can now easily compute the bound 
  $\textup{Var}_{A_2}(d)$ in (\ref{eqnVarB}) for $d=6$ (compare the table in example \ref{exa62}): 
  $$\qquad\qquad \textup{Var}_{A_2}(6) \ = \ \Bigl\lfloor \, \frac{1}{2} \cdot
    \Bigl(\ \underbrace{\frac{15+55-6-14}{2}}_{=10+15} 
    +\underbrace{3 \cdot19 -1 -1}_{=18+19+18}\ \Bigr)\Bigr\rfloor \ = \ 40.\qquad\qquad\Box$$
\end{Exa}

Using some summation formulas we find the following bounds for $d\ge4$.
Some of these are well-known, but we list them because we could not find them
in the literature:
\begin{itemize}
\item $\mu_{A_1}(d) \le \textup{Var}_{A_1}(d) =
  \left\{\begin{array}{ll}
      \frac{23}{48}d^3-\frac{9}{8}d^2+\frac{5}{6}d, & d\equiv 0 \mod 2,\\[0.2em]
      \frac{23}{48}d^{3}-\frac{23}{16}d^{2}+\frac{73}{48}d-\frac{9}{16}, & d\equiv 1 \mod 2.
    \end{array}\right.$\\[0.05em]
\item
  $\mu_{A_2}(d) \le \textup{Var}_{A_2}(d) = \left\{\begin{array}{ll}
      \frac{31}{108}d^{3}-\frac{25}{36}d^{2}+\frac{1}{2}d, & d\equiv 0 \mod 3,\\[0.2em]
      \frac{31}{108}d^{3}-\frac{31}{36}d^{2}+\frac{17}{18}d-\frac{10}{27}, &
      d\equiv 1\mod 3,\\[0.2em]
      \frac{31}{108}d^{3}-\frac{7}{9}d^{2}+\frac{3}{4}d-\frac{5}{27}, &
      d\equiv 2\mod 3.
    \end{array}\right.$\\[0.3em]
\item
  $\mu_{A_3}(d) \le \textup{Var}_{A_3}(d) = \left\{\begin{array}{ll}
      \frac{235}{1152}d^{3}-\frac{49}{96}d^{2}+\frac{13}{36}d, & d\equiv 0
  \mod 4,\\[0.2em] 
      \frac{235}{1152}d^{3}-\frac{235}{384}d^{2}+\frac{785}{1152}d-\frac{35}{128}, & d\equiv 1 \mod 4,\\[0.2em]   
      \frac{235}{1152}d^{3}-\frac{37}{64}d^{2}+\frac{173}{288}d-\frac{3}{16},
      & d\equiv 2\mod 4,\\[0.2em]
     \frac{235}{1152}d^{3}-\frac{209}{384}d^{2}+\frac{569}{1152}d-\frac{35}{384}, &
     d\equiv 3\mod 4.
    \end{array}\right.$\\[0.3em]
\end{itemize} 
The formulas are not correct for $d=3$ for some $j$ because the spectrum of
the $x^3+y^3+z^3=0$ singularity does not have enough spectral numbers to fit into
the description above. 

\subsection{Miyaoka's Bound}
  In \cite{miyp3} Miyaoka gives the following upper bound for the maximum
  number $\mu_{A_j}(d)$ of $A_j$-singularities on a surface of degree $d$ in
  $\dP^3$: 
  \begin{equation}\label{eqnMiyB}
    \mu_{A_j}(d) \le \textup{Miy}_{A_j}(d) := {2\over3}{j+1\over
      j(j+2)}d(d-1)^2 \approx {2\over3}{j+1\over j(j+2)}d^3.
  \end{equation}
  Our variants $\textup{Chm}(M^j_d)$ of Chmutov's
  surfaces give a lower bound of approximately
  $$\mu_{A_j}\bigl(\textup{Chm}(M^j_d)\bigr) \approx {3j+2\over6j(j+1)}d^3$$ such
  singularities for large $d$. 
  This is at least $75\%$ of the best known upper bound:

  \begin{Cor}\label{corQuotient}
    Let $j\in\dN$. 
    For large degree $d$, the quotient of the number of $A_j$-singularities on
    our surfaces $\textup{Chm}(M^j_d)$ and the best known upper bound
    $\textup{Miy}_{A_j}(d)$ is: 
    $$\frac{\mu_{A_j}(\textup{Chm}(M^j_d))}{\textup{Miy}_{A_j}(d)} \approx
    \frac{(j+2)(3j+2)}{4(j+1)^2}.$$ 
    This quotient is greater than $\frac{3}{4}$ for all $j\ge1$, the limit
    for $j\to\infty$ is also $\frac{3}{4}$.\hfill$\Box$
  \end{Cor}

%
%
%
\section{Generalization to Higher Dimensions}
\label{secHighDim}

It is possible to generalize the construction of surfaces with many
$A_j$-singularities described in the previous sections to $\dP^n$, $n\ge4$.  
It turns out that for $n\ge5$, the folding polynomials $F^{A_2}_d(x,y)$ are no
longer the best choice: 
Even for nodal hypersurfaces, the folding polynomials $F^{B_2}_d(x,y)$ lead
to better lower bounds. 

\subsection{Nodal Hypersurfaces in $\dP^n$, $n\ge4$}

As Chmutov mentioned in \cite{chmuP3}, his idea to use the folding
polynomials gives the best lower bounds for the maximum number of nodes on
hypersurfaces in $\dP^4$ of degree $d$ for $d$ large enough.
As Chmutov certainly knew, this can be generalized further to
higher dimensions similar to Givental's construction of cubics in $\dP^n$
\cite[p.\ 419]{AGV2}:
\begin{equation}\label{eqnChmuNewPn}
  \textup{Chm}^n(F_d^{A_2}) : \quad
  \sum_{i=0}^{\lfloor {n-2\over2} \rfloor}(-1)^i F^{A_2}_d(x_{2i},x_{2i+1}) = 
(n \ \textup{mod} \ 2) \smcdot \frac{1}{2}(T_d(x_{n-1}) + 1).
\end{equation}
In some cases, e.g.\ $n=5$, it is better to replace the sign $(-1)^i$ in that
formula by $1$ and to adjust the coefficients on the right-hand side. 
But for $n\ge5$, the asymptotic behaviour (see table
\vref{tabChmPn}) of 
Chmutov's older series $$\textup{TChm}^n_d:\quad \sum_{i=0}^{n-1}
T_d(x_i) = -(n \ \textup{mod} \ 2).$$ 
(see \cite[p.\ 419]{AGV2}) still gives more nodes (exactly
$\bigl(\frac{d-1}{2}\bigr)^n \smcdot {n \choose n/2}$ for odd degree).
The reason for this is that the plane curve $T_d(x)+T_d(y)$ has the critical
values $-1,0,+1$ of which the two non-zero ones sum up to zero.
Nevertheless, for small degree $d$ the hypersurfaces
$\textup{Chm}^n(F_d^{A_2})$ are better. 

In order to improve the asymptotic behaviour of the lower bound slightly, we
can use a folding polynomial associated to another root system. 
Such polynomials were described in \cite{witFoldPoly}, and their critical
points were studied in \cite{sonjaChmuVars} analogous to the case of $A_2$
treated by Chmutov in \cite{chmuP3}.  
It turns out that the folding polynomials $F_d^{B_2}(x,y)$ associated to the root
system $B_2$ are best suited for our purposes. 
They can be defined recursively as follows: $F_0^{B_2} := 1$, \quad  
$F_1^{B_2} := \frac{1}{4}y$, \\[0.3em]
$F_2^{B_2} := \frac{1}{4}y^2-\frac{1}{2}(x^2-2y-4)-1$, \quad 
$F_3^{B_2} := \frac{1}{4}y^3-\frac{3}{4}y(x^2-2y-4) - \frac{3}{4}y,$
\begin{equation}\label{eqnRecFdB2_}F_{d}^{B_2} :=
y(F_{d-1}^{B_2}+F_{d-3}^{B_2}) -
(2+(x^2-2y-4))F_{d-2}^{B_2}-F_{d-4}^{B_2}. 
\end{equation}
These polynomials have exactly three different critical values: 
$-1$, $0$, $+1$.  
The numbers of critical points of $F_{d}^{B_2}$ are: $d\choose2$ with
critical value $0$, $\lfloor \frac{(d-1)}{2}\rfloor \lfloor\frac{d}{2}
\rfloor$ with critical value $-1$. 
The use of these polynomials improves the asymptotic behaviour (for $d$
large) of the best known lower bound for the maximum number of nodes only
slightly.  
In fact, the coefficient of the highest order term does not change (see table
\vref{tabChmPn}). 
Nevertheless, we want to mention: 
\begin{Prop} Let $n\ge2$, $d\ge3$.
Then: $\mu(\textup{Chm}^n(F_d^{B_2})) > \mu(\textup{TChm}^n_d).$
\end{Prop}
\begin{table}
  \begin{center}
    \begin{tabular}{|r|c|c|c|c|c|c|c|c|c|}
      \hline
      \rule{0pt}{1.2em}$n$ & $3$ & $4$ & $5$ & $6$ & $7$ & $8$ & $9$ & $10$\\[0.2em]
      \hline
      \rule{0pt}{1.3em}$\frac{1}{d^n}\smcdot\mu(\textup{Chm}^n(F_d^{A_2})) \approx$ &
      $\frac{5}{12}$ & $\frac{7}{18}$ & $\frac{7}{24}$ & $\frac{19}{72}$ &
      $\frac{35}{144}$ & $\frac{49}{216}$ & $\frac{79}{432}$ & $\frac{25}{144}$\\[0.3em] 
      \hline
      \rule{0pt}{1.4em}$\frac{1}{d^n}\smcdot\mu(\textup{Chm}^n(F_d^{B_2}))
      \approx \frac{1}{d^n}\smcdot\mu(\textup{TChm}^n_d) \approx$ &
      $\frac{3}{2^3}$ & $\frac{3}{2^3}$ & $\frac{5}{2^4}$ & $\frac{5}{2^4}$ &
      $\frac{35}{2^7}$ & $\frac{35}{2^7}$ & $\frac{63}{2^8}$ & $\frac{63}{2^8}$\\[0.3em]   
      \hline
    \end{tabular}\\[0.4em]
    \caption{The asymptotic behaviour of the number of nodes on variants of
      Chmutov's hypersurfaces in $\dP^n$.
      As Chmutov already realized in \cite{chmuP3}, the
      $\textup{Chm}^n(F_d^{A_2})$ are only better for $n=3,4$. 
      For $n\ge5$, the best lower bounds are given by our variant 
      $\textup{Chm}^n(F_d^{B_2})$ which improves Chmutov's oldest examples
      $\textup{TChm}^n_d$ slightly. } 
    \label{tabChmPn}
  \end{center}
\end{table}
It is not true that the folding polynomials $F_d^{A_2}$ and
$F_d^{B_2}$ are the best possible choices in all cases. 
Indeed, for $d=5$, a regular fivegon leads to more nodes. 
For $d=3,4$ there are better constructions for nodal hypersurfaces in $\dP^n$
known \cite{goryQu}. 
In fact, Kalker \cite{kalPhD} already noticed that Varchenko's
upper bound is exact for $d=3$.  

\subsection{Hypersurfaces in $\dP^n$ with $A_j$-Singularities, $j\ge2, n\ge4$} 
\label{subsecHypAj}

Similar to the case of surfaces, we can adapt the equations for the nodal
hypersurfaces to get hypersurfaces $\textup{Chm}^{j,n}(F_d^{B_2})$ (or
$\textup{Chm}^{j,n}(F_d^{A_2})$, $\textup{TChm}^{j,n}_d$) with many
$A_j$-singularities:
\begin{equation}\label{eqnOurPnAj}   
  \textup{Chm}^{j,n}(F_d^{B_2}): \hspace{-0.2cm}
  \sum_{i=0}^{\lfloor {n-3\over2} \rfloor} F^{B_2}_d(x_{2i},x_{2i+1}) = 
  \left\{\begin{array}{@{}l@{\ }l}
      T_d(x_{n-2})+M^j_d(x_{n-1}), & n \, \textup{even}\\
      -\frac{1}{2}(M^j_d(x_{n-1}) + 1), & n \, \textup{odd}.
    \end{array}\right.
\end{equation}
This leads to the asymptotic behaviour given in table \vref{tabChmPnAj}.
Notice that we usually get fewer singularities if we add 
a sign $(-1)^i$ in the sum in contrast to equation (\ref{eqnChmuNewPn}) where
the alternating sign is often better because the folding polynomial
$F_d^{A_2}$ has other critical values than $F_d^{B_2}$.
\begin{table}
  \begin{center}
    \begin{tabular}{|r|c|c|c|c|c|c|c|}
      \hline
      \rule{0pt}{1.2em}$n$ & $3$ & $4$ & $5$ & $6$ & $7$ & $8$ \\[0.2em]
      \hline
      \rule{0pt}{1.4em}$\frac{1}{d^n}\smcdot\mu_{A_2}^n(d) \gtrapprox$ &
      $\frac{2}{9}$ & $\frac{13}{72}$ & $\frac{1}{6}$ & $\frac{13}{96}$ & $\frac{55}{384}$ &
      $\frac{15}{128}$ \\[0.4em]  
      \hline
      \rule{0pt}{1.4em}$\frac{1}{d^n}\smcdot\mu_{A_3}^n(d) \gtrapprox$ & 
      $\frac{11}{72}$ & $\frac{1}{8}$ & $\frac{11}{96}$ & $\frac{3}{32}$ & $\frac{25}{256}$ &
      $\frac{125}{1536}$   \\[0.4em]    
      \hline
      \rule{0pt}{1.4em}$\frac{1}{d^n}\smcdot\mu_{A_4}^n(d) \gtrapprox$ & 
      $\frac{7}{60}$ & $\frac{23}{240}$ & $\frac{7}{80}$ & $\frac{23}{320}$ & $\frac{19}{256}$ &
      $\frac{1}{16}$    \\[0.4em]   
      \hline
      \rule{0pt}{1.4em}$\frac{1}{d^n}\smcdot\mu(\textup{Chm}^{j,n}(F_d^{A_2})) \approx$ & 
      $\frac{3j+2}{6j(j+1)}$ & $\frac{5j+3}{12j(j+1)}$ & $\frac{7j+3}{18j(j+1)}$ &
      $\frac{7j+4}{24j(j+1)}$ & $\frac{19j+16}{72j(j+1)}$ & 
      $\frac{35j+19}{144j(j+1)}$ \\[0.4em]   
      \hline
      \rule{0pt}{1.4em}$\frac{1}{d^n}\smcdot\mu(\textup{Chm}^{j,n}(F_d^{B_2})) \approx$ & 
      $\frac{2j+1}{4j(j+1)}$ & $\frac{3j+2}{8j(j+1)}$ & $\frac{3j+2}{8j(j+1)}$ &
      $\frac{5j+3}{16j(j+1)}$ & $\frac{20j+15}{64j(j+1)}$ & 
      $\frac{35j+20}{128j(j+1)}$ \\[0.4em]   
      \hline
    \end{tabular}\\[0.4em]
    \caption{The asymptotic behaviour of the number of $A_j$-singularities on
      a hypersurface of degree $d$ in $\dP^n$. 
      $\textup{Chm}^{j,n}(F_d^{B_2})$ is better than
      $\textup{Chm}^{j,n}(F_d^{A_2})$ for $n\ge6$.} 
    \label{tabChmPnAj}
  \end{center}
\end{table}

Of course, for small $d,n,j$, it is often easy to write down better lower
bounds. 
E.g., if $n$ is even and $d$ is small, it is often better to replace 
$T_d(x_{n-2})+M^j_d(x_{n-1})$ by a plane curve with the maximum known number
of cusps.
For some specific values of $d$, $j\ge2$, $n\ge4$ there are even better lower
bounds known. 
E.g., Lefschetz \cite{lefCubic5Cusps} constructed a cubic hypersurface in
$\dP^4$ with $5$ cusps which is the maximum possible number. 

%
%
%
\appendix
\section{On Variants of Segre's Construction}

In 1952, B.~Segre \cite{BSegCon2} introduced a construction of surfaces 
with many singularities using pull-back under a branched covering
$\Omega_2^3$. 
Many interesting nodal surfaces can be constructed in this way, e.g., sextics
with $1$ up to $64$ nodes \cite{ccSext} and even Barth's sextic with $65$
nodes \cite{bar65}.  

Shortly after B.~Segre's well-known discovery, Gallarati \cite{galRiflSegre}
generalized this construction to higher dimensions and higher singularities:
\begin{equation}\label{eqnOmegajn}
  \Omega_{j+1}^n: \ \dP^n \to \dP^n, \ (x_0 : x_1 : \cdots : x_n) \mapsto
  (x_0^{j+1} : x_1^{j+1} : \cdots : x_n^{j+1}), \quad j\in\dN.
\end{equation}
Gallarati does not give a general formula for the number and type of
singularities one obtains using this map.
He only computes some examples. 
But it is easy to derive a formula for hypersurfaces with
$A_{j}$-singularities similar to B.~Segre's case of nodal surfaces in
$\dP^3$:
Let $F_0$ be a hypersurface in $\dP^n$ of degree $d_0$ with $k_0$ singularities
of type $A_{j}$. 
Take $n+1$ general hyperplanes tangent to $F_0$ as the coordinate
$(n+1)$-hedron.  
The degree of the map $\Omega_{j+1}^n$ is $(j+1)^n$ away from the coordinate
hyperplanes.
It is $(j+1)^{n-1}$ on a general intersection point of two of the coordinate
hyperplanes, and $(j+1)^{n-i}, i=2,3,\dots,n,$ for even more special points on the
coordinate hyperplanes. 
For our generic choice of coordinate hyperplanes tangent to $F_0$
the pull-back under $\Omega_{j+1}^n$ thus gives a hypersurface $F_1$ 
in $\dP^n$ of degree $d_1 := (j+1)\smcdot d_0$ with 
\begin{equation}\label{eqnF_1}
  \mu_{A_{j}}(F_1) = (j+1)^n\smcdot k_0 + (n+1)\smcdot (j+1)^{n-1}
\end{equation}
singularities of type $A_{j}$. 
Applying the same construction to $F_1$, we obtain a hypersurface $F_2$ in
$\dP^n$ of degree $d_2 := (j+1)^2\smcdot d_0$ with 
$$\mu_{A_{j}}(F_2) = (j+1)^n\bigl((j+1)^n\smcdot k_0 + (n+1)\smcdot (j+1)^{n-1}\bigr) +
  (n+1)(j+1)^{n-1}$$
singularities of type $A_{j}$. 
Iterating this, we get a hypersurface $F_i$ of
  degree $d_i := (j+1)^i\smcdot d_0$ with 
$$\mu_{A_{j}}(F_i) = (j+1)^{ni}\smcdot k_0 +
 \frac{n+1}{j+1}\cdot\Bigl(\frac{(j+1)^{n(i+1)}-1}{(j+1)^n-1} - 1\Bigr)$$
  singularities of type $A_{j}$.  
Asymptotically, we thus have:
\begin{equation}\label{eqnApproxF_i}
\mu_{A_{j}}(F_i)
\approx \frac{1}{d_0^n}\smcdot\left(k_0 +
  \frac{(n+1)\smcdot(j+1)^{n-1}}{(j+1)^n-1}\right)\smcdot d_i^n \quad
\textup{for} \ i \ \textup{large}. 
\end{equation}
For $n\ge3$, this lower bound is asymptotically not as good as ours 
presented in the main text.
But for low degree, B.~Segre's method sometimes gives more singularities:
E.g., when applying  
$\Omega_{j+1}^3$ to a smooth quadric, (\ref{eqnF_1}) yields to: 
\begin{Cor}\label{ThmSext36}
  Let $j\in\dN$. 
  There exist surfaces of degree $d=2\smcdot(j+1)$ with $4\smcdot(j+1)^2$
  singularities of type $A_j$.
\end{Cor} 

E.g., for $n=3, j=2$, we obtain $\mu_{A_2}(6)\ge36$, and with Miyaoka's upper
bound: $36\le\mu_{A_2}(6)\le37$. 
For $n=2$, our construction presented in subsection \vref{subsecHypAj} only
leads to plane curves of degree $d$ with $\approx \frac{1}{4}\smcdot d^2$
cusps whereas the generalization of B.~Segre's construction gives $\approx
\frac{9}{32}\smcdot d^2$ such singularities when starting with a smooth
conic. 
This idea was taken up later by several people. 
To our knowledge, the currently best result is due to Vik.~S. Kulikov 
\cite{kulGenChis}.
He used a quartic with three cusps as a starting point. 
At every other iteration step he was able to choose a bitangent to the curve as
one of the coordinate axes. 
This yields to approximately $\frac{283}{60\cdot16}\smcdot d^2$ cusps. 
So, in the case of plane curves of degree $d$, variants of B.~Segre's idea
still give the best known general lower 
bound for the maximum number $\mu_{A_2}^2(d)$ of cusps.

In higher dimensions, our construction gives a better lower bound than
this generalization of B.~Segre's construction.
Notice that it might be able to adapt B.~Segre's construction similar to
the case of curves: 
In the case of surfaces, it might be possible to choose triple tangent planes
as coordinate planes. 
But even when starting from a $36$-cuspidal sextic this
would yield to surfaces with less cusps.  

Finally, we want to mention that it is easy to compute how many singularities
we need to improve the best known lower bounds using the formula
(\ref{eqnApproxF_i}).   
Let us look at nodal surfaces: To improve Chmutov's lower bound $\approx
\frac{5}{12}d^3$ for the maximum number of nodes on a surface of degree $d$,
it suffices to construct a surface of degree $d_0$ with $k_0$ nodes, s.t.\
$k_0 > \frac{5}{12} d_0^3 - \frac{16}{7}$. 
Comparing this with Miyaoka's upper bound, we find, e.g., that a $13652$-nodal
surface of degree $32$ or a $109225$-nodal surface of degree $64$ would be
sufficient. 

\bibliographystyle{plain} 
\bibliography{papers}

\end{document}